\documentclass{math}
\usepackage{hyperref}
\hypersetup{
colorlinks=true,
urlcolor=blue,
citecolor=blue}
\usepackage[all]{xy}
\usepackage{amsfonts,amssymb,amsmath,amsgen,amsopn,amsbsy,theorem,graphicx,epsfig}
\usepackage{amscd,bezier,latexsym,mathrsfs,enumerate}\usepackage[utf8]{inputenc}\usepackage[english]{babel}
\usepackage[dvipsnames]{xcolor}

\year{2018}
\vol{}
\fpage{}
\lpage{}
\doi{}

\title{A hybrid simulation for singularly perturbed system of two-point reaction-diffusion boundary-value problems}
\author[CENGIZCI et al.]{
\textbf{S\"uleyman CENGIZCI$^{1}$\thanks{Correspondence: suleyman.cengizci@antalya.edu.tr}~, Srinivasan NATESAN$^{2}$, Mehmet Tarık ATAY$^{3}$}\\
$^{1}$Department of Computer Programming, Antalya Bilim University, Antalya, Turkey\\
$^{2}$Department of Mathematics, Indian Institute of Technology Guwahati, Guwahati, India \\
$^{3}$Department of Mechanical Engineering, Abdullah G\"ul University, Kayseri, Turkey\\

\\ [1.8em]

\rec{28.07.2018}
\acc{.201}
\finv{..201}
}

\amssayisi{2010 {\itshape AMS Mathematics Subject Classification:} 65L10, 65L11 }

\newcommand{\bc}{\begin{center}}
\newcommand{\ec}{\end{center}}

\numberwithin{equation}{section}

\newtheorem{theorem}{Theorem}[section]

\newtheorem{corollary}[theorem]{Corollary}

\newtheorem{example}[theorem]{Example}

\newtheorem{lemma}[theorem]{Lemma}

\begin{document}

\maketitle
\begin{abstract}
This study concerns with singularly perturbed systems of second-order
reaction-diffusion equations in ODE's. To handle this type of problems, a
numerical-asymptotic hybrid method is employed. In this hybrid method, an
efficient asymptotic method, the so-called Successive complementary expansion method (SCEM) is applied first and then, a numerical method based on finite differences is proposed to approximate the solution of the corresponding singularly perturbed reaction-diffusion systems. Two illustrative examples are provided to show the efficiency and easy-applicability of the present method with convergence properties.
\end{abstract}

\section{Introduction}\label{Sec:1}

System of differential equations arise in many branches of science such as modelling of electrical networks, mechanical systems, marketing problems,
earthquake/tsunami problems, chemotaxis processes, semiconductor physics,
etc. Because of the physical importance, many studies were devoted to this area. In \cite{11}, system of differential equations are studied in detail and in \cite{2,4} some numerical treatments are examined. On the other hand, singularly perturbed differential equations that involve positive small perturbation parameter(s) $0<\epsilon \ll 1$ multiplied with the highest order derivative term are also important concept of applied sciences and it is a well-known fact that standard numerical techniques are often insufficient to handle them. Some application areas may be given as control theory, fluid mechanics, quantum mechanics, combustion theory, signal and image processing, pharmacokinetics, etc. One can find theoretical considerations on singular perturbation problems in \cite{14,22,24,28} and various approximation methods in \cite{5,8,12,13,15,16,29}.

The present paper concerns with the approximations of the solution of
singularly perturbed systems of reaction-diffusion boundary-value problems (BVPs) that frequently arise in electroanalytical chemistry and population dynamics problems.\ In recent years various methods were employed to obtain approximations to the solution of this kind of problems. In \cite{3,19}, finite difference methods (FDM) for singularly perturbed convection-diffusion systems, in \cite{17,18} finite element method (FEM)
and in \cite{10,20,23,26} finite difference methods for singularly
perturbed reaction-diffusion systems are examined. For singularly perturbed system of reaction-diffusion BVPs, in \cite{200}, the authors provided a robust computational technique , whereas in \cite{110}, optimal order error estimates are obtained on equidistributed grids.

In this paper, we employ a numerical-asymptotic hybrid method for approximating the solution of singularly perturbed system of two-point boundary-value problems reaction-diffusion type. At the first step, an efficient asymptotic method that is introduced in \cite{21}, called Successive Complementary Expansion Method (SCEM), and later a finite difference method given in \cite{25} are applied.

In order to examine the convergence of approximations to the solution of
singularly perturbed system of differential equations, we need to define the norm that we will use in the remainder parts of this study. In \cite{9} one can find that the appropriate norm for this examination is the \textit{maximum norm} that is given by
\begin{equation*}
\Vert r\parallel =\underset{\overline{\Omega }}{\max }\left\vert r\left(
x\right) \right\vert ,\text{ }\Vert \overrightarrow{r}\parallel =\underset{%
\overline{\Omega }}{\max }\left\{ \Vert r_{i}\parallel \right\} ,
\end{equation*}%
where $\overrightarrow{r}=\left( r_{1},r_{2},\cdots,r_{n}\right) .$

The remainder sections of this paper continue as follows: In Section \ref{sec2}, the continuous problem and its general properties are explained. In section \ref{sec3}, the numerical-asymptotic hybrid method that we employ to solve system of singularly perturbed two-point reaction-diffusion equations is described. In section \ref{sec4}, two illustrative examples are provided to show the implementation of the method and in the last section, we discuss the findings.

\goodbreak
\section{The continuous problem}\label{sec2}

In general, a coupled system of singularly perturbed reaction-diffusion
ODE's is given as
\begin{equation}\label{sys1}
\left\{
\begin{array}{ll}
-\varepsilon y_{1}^{\prime \prime }\left( x\right) +a_{11}\left( x\right)
y_{1}\left( x\right) +a_{12}\left( x\right) y_{2}\left( x\right)
 = f_{1}\left( x\right) , \\[6pt]
-\eta y_{2}^{\prime \prime }\left( x\right) +a_{21}\left( x\right)
y_{1}\left( x\right) +a_{22}\left( x\right) y_{2}\left( x\right)
= f_{2}\left( x\right) ,
\end{array}\right.
\end{equation}
with some suitable boundary or initial conditions. The presence of positive small parameters $0<\varepsilon \ll 1$ and $0<\eta \ll 1$ causes rapid changes (boundary layers) near the end-points of the domain $\overline{\Omega }.$ These layer behaviors can be examined in three different cases:

\goodbreak\noindent
\textbf{Case 1:}\textit{\ }$\mathbf{0<\varepsilon <\eta \ll 1}$ In this case, both components of the solution have boundary layers of width $O\left( \eta \ln \eta \right) $ and $y_{1}\left( x\right) $ has an additional sublayer of width $O\left( \varepsilon \ln \varepsilon \right)$.

\goodbreak\noindent
\textbf{Case 2:}\textit{\ }$\mathbf{0<\varepsilon \ll 1}$ and $\mathbf{\eta =1}$ In this case, only the first component of the solution has boundary layers of width $O\left( \varepsilon \ln \varepsilon \right)$.

\goodbreak\noindent
\textbf{Case 3:}\textit{\ }$\mathbf{0<\varepsilon =\eta \ll 1}$ In this case, both components of the solution have boundary layers of width $O\left( \varepsilon \ln \varepsilon \right) $ (or equivalently $O\left( \eta \ln \eta\right))$.

We are interested in Case 2 - Case 3 and in finding $\overrightarrow{y}%
\left( x\right) \in C^{2}\left( \overline{\Omega }\right) ,$ such that for all $%
x\in \Omega =\left( 0,1\right) $ and for the problem
\begin{equation}\label{sys2}
\left\{
\begin{array}{ll}
-\varepsilon y_{1}^{\prime \prime }\left( x\right) +a_{11}\left( x\right)
y_{1}\left( x\right) +a_{12}\left( x\right) y_{2}\left( x\right)
 = f_{1}\left( x\right) , \\[6pt]
-\varepsilon y_{2}^{\prime \prime }\left( x\right) +a_{21}\left( x\right)
y_{1}\left( x\right) +a_{22}\left( x\right) y_{2}\left( x\right)
 = f_{2}\left( x\right) ,   \\[6pt]
y_{1}\left( 0\right) =y_{1}\left(
1\right) =0, \quad y_{2}\left( 0\right) =y_{2}\left( 1\right) =0.
\end{array}\right.
\end{equation}

The system of BVPs given in (\ref{sys2}) can be rewritten in the matrix form as
\begin{equation}\label{mat1}
\left\{
\begin{array}{ll}
L_{\varepsilon }\overrightarrow{y}\left( x\right) =%
\begin{bmatrix}
-\varepsilon \frac{d^{2}}{dx^{2}} & 0 \\
0 & -\varepsilon \frac{d^{2}}{dx^{2}}%
\end{bmatrix}%
\overrightarrow{y}\left( x\right) +A\left( x\right) \overrightarrow{y}\left(
x\right) =\overrightarrow{f}\left( x\right), \\[10pt]
\overrightarrow{y}\left( 0\right) =\left[
y_{1}\left( 0\right) ,y_{2}\left( 0\right) \right] ^{T}, \quad
\overrightarrow{y}%
\left( 1\right) =\left[ y_{1}\left( 1\right) ,y_{2}\left( 1\right) \right]^{T}
\end{array}\right.
\end{equation}
where $A\left( x\right) =%
\begin{bmatrix}
a_{11}\left( x\right) & a_{12}\left( x\right) \\
a_{21}\left( x\right) & a_{22}\left( x\right)%
\end{bmatrix}%
$ and $\overrightarrow{f}\left( x\right) =%
\begin{bmatrix}
f_{1}\left( x\right) \\
f_{2}\left( x\right)%
\end{bmatrix}%
.$ The functions $a_{ij}\left( x\right) ,$ $f_{i}\left( x\right) \in
C^{2}\left( \Omega \right) $ for $i,j=1,2.$ Now, we should impose two
assumptions on problem (\ref{mat1}). First one is for strictly
diagonally dominance of matrix $A$ and the second one is to be able to
establish maximum principle theorem:

\textbf{Assumption 1}
$a_{11}\left( x\right) >\left\vert a_{12}\left( x\right)
\right\vert $ and $a_{22}\left( x\right) >\left\vert a_{21}\left( x\right)
\right\vert ,$ for all $x\in \overline{\Omega },$

\textbf{Assumption 2}
$a_{12}\left( x\right) \leq 0$ and $a_{21}\left(
x\right) \leq 0,$ for all $x\in \overline{\Omega }.$

The proof of the following lemma and the corollary can be seen in the paper \cite{20}.

\begin{lemma}(Maximum Principle)
Consider the system of singularly
perturbed BVP's (\ref{mat1}). If $\overrightarrow{y}\left( 0\right) \geqslant
\overrightarrow{0},$ $\overrightarrow{y}\left( 1\right) \geqslant
\overrightarrow{0}$ and $L_{\varepsilon }\overrightarrow{y}\left( x\right)
\geqslant \overrightarrow{0}$ for all $x\in $ $\Omega ,$ then $\overrightarrow{y%
}\left( x\right) \geqslant \overrightarrow{0}$,  $\forall x\in
\overline{\Omega }.$
\end{lemma}

\begin{corollary}(Stability)
If $\overrightarrow{y}\left( x\right) $ is
the solution of (2.3), then the stability bound inequality
\begin{equation*}
\Vert \overrightarrow{y}\left( x\right) \Vert \leq \frac{1}{\delta }\Vert
\overrightarrow{f}\Vert +\Vert \overrightarrow{y}\left( 0\right) \Vert
+\Vert \overrightarrow{y}\left( 1\right) \Vert
\end{equation*}%
holds, where $\delta =\underset{\overline{\Omega }}{\min }\left\{ a_{11}\left(
x\right) +a_{12}\left( x\right) ,a_{21}\left( x\right) +a_{22}\left(
x\right) \right\} .$
\end{corollary}

Under these above-mentioned assumptions and conditions, the hybrid method can be given as in the next section.

\section{The hybrid method}\label{sec3}

In this section, we first give a short overview of asymptotic approximations and then explain the hybrid method by which we obtain highly accurate approximations to the systems of singularly perturbed differential equations.

Let $E$ be a set of all real-valued functions that depend on $\varepsilon ,$ strictly positive and continuous in $(0,\varepsilon _{0}]$ and such that $\displaystyle\lim_{\varepsilon \rightarrow 0}\delta (\varepsilon )$ exists and for each $\delta _{1}\left( \varepsilon \right) $,$\delta _{2}\left( \varepsilon \right) \in E$, $\delta _{1}\left( \varepsilon \right) \delta _{2}\left( \varepsilon \right) \in E$ holds. A function $\delta _{i}(\varepsilon )$ that satisfies these conditions are called \textit{order function}. Given two functions $\phi (x,\varepsilon )$ and $\phi _{a}(x,\varepsilon )$ defined in a domain $\overline{\Omega }$ are asymptotically identical to order $\delta \left( \varepsilon \right)
$ if their difference is asymptotically smaller than $\delta \left( \varepsilon \right) $, where $\delta \left( \varepsilon \right) $ is an order function, that is,
\begin{equation}
\phi (x,\varepsilon )-\phi _{a}(x,\varepsilon )=o(\delta (\varepsilon ))
\end{equation}
where $\varepsilon $ is the small parameter arising from the physical problem under consideration. The function $\phi _{a}(x,\varepsilon )$ is named as \textit{asymptotic approximation} of the function $\phi (x,\varepsilon ).$ Asymptotic approximations in general form defined by
\begin{equation}\label{asymp1}
\phi _{a}(x,\varepsilon )=\sum_{i=1}^{n}\delta _{i}(\varepsilon )\varphi
_{i}(x,\varepsilon )
\end{equation}
where the asymptotic sequence of order functions $\delta _{i}(\varepsilon )$ satisfy the condition $\delta _{i+1}(\varepsilon )=o(\delta _{i}(\varepsilon ))$, as $\varepsilon \rightarrow 0$. Under these conditions, the approximation (\ref{asymp1}) is named as \textit{generalized asymptotic expansion}. If the expansion (\ref{asymp1}) is written in the form of%
\begin{equation}\label{asymp2}
\phi _{a}(x,\varepsilon )=E_{0}\phi =\sum_{i=1}^{n}\delta _{i}^{(0)}(\varepsilon
)\varphi _{i}^{(0)}(x),
\end{equation}
then it is called as \textit{regular asymptotic expansion} where the special operator $E_{0}$ is \textit{outer\ expansion operator} at a given order $\delta (\varepsilon ).$ Thus $\phi -E_{0}\phi =o(\delta (\varepsilon ))$. Interesting cases occur when the function is not regular in $\overline{\Omega }$ so (\ref{asymp1}) or (\ref{asymp2}) is valid only in a restricted region $\overline{\Omega }_{0} \in $ $\overline{\Omega }$ called the outer region. This is a singular perturbation problem and we must introduce boundary layer domains. We introduce an inner domain which can be formally denoted by $\overline{\Omega }_{1}= \overline{\Omega }-\overline{\Omega }_{0}$ and corresponding inner layer variable
located near the point $x=x_{0}$ as $\overline{x}=\frac{x-x_{0}}{\xi
(\varepsilon )}$, $\xi (\varepsilon )$ being the order of thickness of this boundary layer. If a regular expansion can be constructed in $\overline{\Omega }_{1}$, we can write down
\begin{equation}\label{asymp3}
\phi _{a}(x,\varepsilon )=E_{1}\phi =\sum_{i=1}^{n}\delta _{i}^{(1)}(\varepsilon
)\varphi _{i}^{(1)}(\overline{x})
\end{equation}
where the \textit{inner expansion operator} $E_{1}$ is defined in $\overline{\Omega }_{1}$ at the same order of $\delta (\varepsilon )$ as the outer expansion operator $E_{0};$ thus, $\phi -E_{1}\phi =o(\delta (\varepsilon ))$ and
\begin{equation*}\label{asymp4}
\phi _{a}=E_{0}\phi +E_{1}\phi -E_{1}E_{0}\phi
\end{equation*}%
is clearly uniformly valid approximation (UVA) \cite{1,6,7}. For only
one singularly perturbed differential equation, the uniformly valid SCEM
approximation is in the regular form given by
\begin{equation}\label{asymp5}
y_{n}^{scem}(x,\overline{x},\varepsilon )=\overset{n}{\underset{i=1}{\sum }}%
\delta _{i}(\varepsilon )\left[ y_{i}(x)+\Psi _{i}(\overline{x})\right]
\end{equation}
where$\ \left\{ \delta _{i}\left( \varepsilon \right) \right\} \ $is an
asymptotic sequence and functions $\Psi _{i}(\overline{x})$ are the
complementary functions that depend on $\overline{x}$. If the functions $
y_{i}(x)$\ and $\Psi _{i}(\overline{x})$ depend also on $\varepsilon $, the uniformly valid SCEM approximation is called as \textit{generalized SCEM approximation }and given by
\begin{equation}\label{scem}
y_{ng}^{scem}(x,\overline{x},\varepsilon )=\overset{n}{\underset{i=1}{\sum }} \delta _{i}(\varepsilon )\left[ y_{i}(x,\varepsilon )+\Psi _{i}(\overline{x},\varepsilon )\right] .
\end{equation}
Since the problem (\ref{mat1}) exhibits boundary layers at both the end-points of the interval $\Omega = (0,1)$, in general, the generalized SCEM approximation will be adopted as follows:

For the outer region, that is far from the end-points, the outer approximation will be in the form of
\begin{equation}\label{outer1}
\overrightarrow{y}^{out}\left( x,\delta \right) =\left[
\begin{array}{c}
\overrightarrow{y}_{1}^{out}\left( x,\delta \left( \varepsilon \right) \right) \\[6pt]
\overrightarrow{y}_{2}^{out}\left( x,\delta \left( \eta \right) \right)%
\end{array}
\right] =\left[
\begin{array}{c}
\overrightarrow{y}_{1}^{out\left( 1\right) }\left( x\right) +\delta \left(
\varepsilon \right) \overrightarrow{y}_{1}^{out\left( 2\right) }\left( x\right)
+\delta ^{2}\left( \varepsilon \right) \overrightarrow{y}_{1}^{out\left(
3\right) }\left( x\right) +\cdots \\
\overrightarrow{y}_{2}^{out\left( 2\right) }\left( x\right) +\delta \left(
\eta \right) \overrightarrow{y}_{2}^{out\left( 2\right) }\left( x\right)
+\delta ^{2}\left( \eta \right) \overrightarrow{y}_{2}^{out\left( 3\right)
}\left( x\right) +\cdots%
\end{array}
\right].
\end{equation}
If equation (\ref{outer1}) is substituted into (\ref{mat1}) and the powers of $\varepsilon $ and $\eta $ are balanced, one gets the asymptotic approximation for outer region. In order to cope with the approximational difficulties at the end-points, stretching variables will be introduced as $\overline{x}_{L}= \frac{x}{\sqrt{\varepsilon }}$ for the left-end and $\overline{x}_{R}=\frac{x-1}{\sqrt{\varepsilon }}$ for the right-end. Using these transformations with the help of chain rule and substituting into the equation (\ref{mat1}) one obtains the complementary functions as the solution of following sub-problems
\begin{equation}\label{subprob1}
\left\{
\begin{array}{ll}
-\left(\Psi _{1}^{Lcomp}\right)'' \left( x,\overline{x}_{L},\varepsilon \right) + a_{11}\left( \overline{x}_{L}\right) \Psi _{1}^{Lcomp}\left( x,\overline{x}_{L},\varepsilon \right) +a_{12}\left( \overline{x}_{L}\right) \Psi
_{2}^{Lcomp}\left( x,\overline{x}_{L},\varepsilon \right)  = f_{1}\left(
\overline{x}_{L}\right) , \\[8pt]
-\left(\Psi_{2}^{Lcomp}\right)''\left( x,\overline{x}_{L},\varepsilon \right) + a_{21}\left( \overline{x}_{L}\right) \Psi _{1}^{Lcomp}\left( x,\overline{x}_{L},\varepsilon \right) +a_{22}\left( \overline{x}_{L}\right) \Psi_{2}^{Lcomp}\left( x,\overline{x}_{L},\varepsilon \right) = f_{2}\left(
\overline{x}_{L}\right),\\[8pt]
\Psi _{1}^{Lcomp}\left( 0,0,\varepsilon \right)
=-y_{1}^{out(}\left( 0\right), \quad \Psi _{1}^{Lcomp}\left( 1,\frac{1}{\sqrt{\varepsilon }},\varepsilon \right) =-y_{1}^{out}\left( 1\right),\\[8pt]
\Psi_{2}^{Lcomp}\left( 0,0,\varepsilon \right) =-y_{2}^{out(}\left( 0\right), \quad
\Psi _{2}^{Lcomp}\left( 1,\frac{1}{\sqrt{\varepsilon }},\varepsilon \right) = -y_{2}^{out}\left( 1\right)
\end{array}\right.
\end{equation}
and
\begin{equation}\label{subprob2}
\left\{
\begin{array}{ll}
-\left(\Psi_{1}^{Rcomp}\right)''\left( x,\overline{x}_{R},\varepsilon \right) + a_{11}\left( \overline{x}_{R}\right) \Psi _{1}^{Rcomp}\left( x,\overline{x}_{R},\varepsilon \right) +a_{12}\left( \overline{x}_{R}\right) \Psi_{2}^{Rcomp}\left( x,\overline{x}_{R},\varepsilon \right) = f_{1}\left(
\overline{x}_{R}\right) , \\[8pt]
-\left(\Psi_{2}^{Rcomp}\right)''\left( x,\overline{x}_{R},\varepsilon \right) + a_{21}\left( \overline{x}_{R}\right) \Psi _{1}^{Rcomp}\left( x,\overline{x}_{R},\varepsilon \right) +a_{22}\left( \overline{x}_{R}\right) \Psi_{2}^{Rcomp}\left( x,\overline{x}_{R},\varepsilon \right) = f_{2}\left(
\overline{x}_{R}\right), \\[8pt]
\Psi _{1}^{Rcomp}\left( 0,\frac{-1}{\sqrt{\varepsilon }},\varepsilon \right) =-y_{1}^{out}\left( 0\right), \quad \Psi_{1}^{Rcomp}\left( 1,0,\varepsilon \right) =-y_{1}^{out}\left( 1\right), \\[8pt]
\Psi _{2}^{Rcomp}\left( 0,\frac{-1}{\sqrt{\varepsilon }},\varepsilon \right) = -y_{2}^{out}\left( 0\right), \quad
\Psi_{1}^{Rcomp}\left( 1,0,\varepsilon\right) =-y_{2}^{out}\left( 1\right).
\end{array}\right.
\end{equation}
where the superscripts $Lcomp(i)$ and $Rcomp(i)$ denote the $i^{th}$ complementary approximations of left and right layer problems respectively, and subscripts $1,2$ denote the first and second components of the approximations. If an asymptotic approximation for the complementary functions in the following form is adopted
\begin{equation*}
\Psi \left( x,\overline{x},\varepsilon \right) =\left[
\begin{array}{c}
\Psi _{1}\left( x,\overline{x},\delta \left( \varepsilon \right) \right) \\
\Psi _{2}\left( x,\overline{x},\delta \left( \eta \right) \right)%
\end{array}%
\right] =\left[
\begin{array}{c}
\Psi _{1}^{\left( 1\right) }\left( \overline{x},\delta \left( \varepsilon
\right) \right) +\delta \left( \varepsilon \right) \Psi _{1}^{\left( 2\right)}\left( \overline{x},\delta \left( \varepsilon \right) \right) +\delta^{2}\left( \varepsilon \right) \Psi _{1}^{\left( 3\right) }\left( \overline{x},\delta \left( \varepsilon \right) \right) +\cdots \\
\Psi _{2}^{\left( 1\right) }\left( \overline{x},\delta \left( \varepsilon
\right) \right) +\delta \left( \eta \right) \Psi _{2}^{\left( 2\right)
}\left( \overline{x},\delta \left( \varepsilon \right) \right) +\delta
^{2}\left( \eta \right) \Psi _{2}^{\left( 3\right) }\left( \overline{x}%
,\delta \left( \varepsilon \right) \right) +\cdots%
\end{array}
\right] ,
\end{equation*}
and substituted to problem (\ref{mat1}) one gets the asymptotic approximations for complementary functions again balancing them with respect to the power of $\varepsilon $ and $\eta$.

To this end, the first iteration of the hybrid method is in the form of
\begin{equation}
\overrightarrow{y}^{(1)}\left( x,\overline{x}_{L},\overline{x}_{R},\varepsilon
\right) =\overrightarrow{y}^{out(1)}\left( x,\varepsilon \right) +\left[ \frac{%
\overrightarrow{\Psi }^{Lcomp(1)}\left( x,\overline{x}_{L},\varepsilon \right) +%
\overrightarrow{\Psi }^{Rcomp(1)}\left( x,\overline{x}_{R},\varepsilon \right)
}{2}\right] ,
\end{equation}
where complementary functions are numerically solved by the numerical method that is given in \cite{25} based on finite differences and that implements the three-stage Lobatto IIIa formula.

\goodbreak
\section{Illustrative Examples}\label{sec4}

In this section, two numerical examples will be solved. In the first one, all the process explained in detail and, all the computations are performed in {\sc Matlab2015b} using double precision.

\begin{example}\label{exam1}
Consider the system of coupled singularly perturbed
reaction-diffusion equations \cite{27}:
\begin{equation}\label{exam_pb1}
\left\{
\begin{array}{ll}
-\varepsilon y_{1}^{\prime \prime }\left( x\right) +4y_{1}\left( x\right)
-2y_{2}\left( x\right) = 1, \quad x \in \Omega = (0,1)\\[6pt]
-\varepsilon y_{2}^{\prime \prime }\left( x\right) -y_{1}\left( x\right)
+3y_{2}\left( x\right) = 2,\\[6pt]
y_{1}\left( 0\right) =y_{1}\left( 1\right) =0, \quad
y_{2}\left( 0\right) =y_{2}\left( 1\right) =0.
\end{array}\right.
\end{equation}
\end{example}
As one can see in Figure 1, solution of this problem exhibits boundary layer behavior at the end-points of the interval as, $\varepsilon \rightarrow 0.$ Therefore, for both end-point a stretching variable will be introduced as $\overline{x}_{L}=%
\frac{x}{\sqrt{\varepsilon }}$ for the left-end and $\overline{x}_{R}=\frac{x-1%
}{\sqrt{\varepsilon }}$ for the right-end. But first, the reduced problem
should be obtained taking $\varepsilon =0$:
\begin{equation}\label{outer1}
\left\{
\begin{array}{ll}
4y_{1}^{out(1)}\left( x,\varepsilon \right) -2y_{2}^{out(1)}\left( x,\varepsilon\right) = 1, \\[8pt]
-y_{1}^{out(1)}\left( x,\varepsilon \right) +3y_{2}^{out(1)}\left( x,\varepsilon \right) = 2,
\end{array}\right.
\end{equation}
and it is obvious that the solution to this reduced system is $
y_{1}^{out(1)}\left( x\right) =0.7$ and $y_{2}^{out(1)}\left( x\right) =0.9,$ where the superscript $out(i)$ denotes the $i^{th}$ approximation to the outer layer problem. For the left and right inner layer problems, adopting the stretching variables $\overline{x}_{L}$ and $\overline{x}_{R}$ respectively, one gets the systems
\begin{equation}\label{first1}
\left\{
\begin{array}{ll}
-\left(\Psi_{1}^{Lcomp(1)}\right)''\left( x,\overline{x}_{L},\varepsilon
\right) +4\Psi _{1}^{Lcomp(1)}\left( x,\overline{x}_{L},\varepsilon \right)-2\Psi _{2}^{Lcomp(1)}\left( x,\overline{x}_{L},\varepsilon \right) = 1, \\[8pt]
-\left(\Psi_{2}^{Lcomp(1)}\right)''\left( x,\overline{x}_{L},\varepsilon
\right) -\Psi _{1}^{Lcomp(1)}\left( x,\overline{x}_{L},\varepsilon \right) + 3\Psi_{2}^{Lcomp(1)}\left( x,\overline{x}_{L},\varepsilon \right) = 2, \\[8pt]
\Psi _{1}^{Lcomp(1)}\left( 0,0,\varepsilon\right) =-0.7, \quad
\Psi _{1}^{Lcomp(1)}\left( 1,\frac{1}{\sqrt{\varepsilon }},\varepsilon \right) =-0.7, \\[8pt]
\Psi _{2}^{Lcomp(1)}\left( 0,0,\varepsilon \right)=-0.9, \quad
\Psi _{2}^{Lcomp(1)}\left( 1,\frac{1}{\sqrt{\varepsilon }},\varepsilon\right) =-0.9
\end{array}\right.
\end{equation}
and
\begin{equation}\label{first2}
\left\{
\begin{array}{ll}
-\left(\Psi_{1}^{Rcomp(1)}\right)''\left( x,\overline{x}_{R},\varepsilon
\right) +4\Psi _{1}^{Rcomp(1)}\left( x,\overline{x}_{R},\varepsilon \right) - 2\Psi _{2}^{Rcomp(1)}\left( x,\overline{x}_{R},\varepsilon \right) = 1, \\[8pt]
-\left(\Psi_{2}^{Rcomp(1)}\right)'' \left( x,\overline{x}_{R},\varepsilon \right) -\Psi _{1}^{Rcomp(1)}\left( x,\overline{x}_{R},\varepsilon \right) + 3\Psi _{2}^{Rcomp(1)}\left( x,\overline{x}_{R},\varepsilon \right) = 2, \\[8pt]
\Psi _{1}^{Rcomp(1)}\left( 0,\frac{-1}{\sqrt{\varepsilon }},\varepsilon \right) =-0.7, \quad \Psi_{1}^{Rcomp(1)}\left(1,0,\varepsilon \right) =-0.7, \\[8pt]
\Psi _{2}^{Rcomp(1)}\left( 0,\frac{-1}{\sqrt{\varepsilon }},\varepsilon \right) =-0.9, \quad \Psi _{2}^{Rcomp(1)}\left(1,0,\varepsilon \right) =-0.9
\end{array}\right.
\end{equation}
In order to employ the numerical method \cite{25} and code the problem in
{\sc Matlab}, we should transform these new systems into first order systems. The transformed problem corresponding to problem (\ref{first1}) can be coded using {\sc Matlab bvp4c} as
\begin{verbatim}
function res = twobc(ya,yb)
res = [ ya(1)+7/10
yb(1)+7/10
ya(3)+9/10
yb(3)+9/10];
function dydx = twoode(x,y)
dydx = [y(2)-1+4*y(1)-2*y(3) y(4)-2-y(1)+3*y(3)];
ep=0.01;
solinit = bvpinit(linspace(0,1/sqrt(ep),1000),[1 1 1 1]);
sol = bvp4c(@twoode,@twobc,solinit);
x=0:0.01:1/sqrt(ep);
y = deval(sol,x);
double(y);
y11=0.5*(y(1,:)+0.7);
y33=0.5*(y(3,:)+0.9);
\end{verbatim}

The double mesh principle will be used for estimating the maximum point-wise errors and computing the rate of convergence in the computed approximations.
\[
D_{\varepsilon ,i}^{N}=\underset{x_{j}\in \overline{\Omega }_{\varepsilon }^{N}}{\max}\mid Y_{i}^{2N}\left( x_{j}\right) -Y_{i}^{N}\left( x_{j}\right) \mid
\]
for $i=1,2$ and $D_{i}^{N}=\underset{\varepsilon }{\max }D_{\varepsilon ,i}^{N}$, where $Y_{i}^{N}\left( x_{j}\right) $ and $Y_{i}^{2N}\left( x_{j}\right) $ denote the computed approximations at the point $x=x_{j},$ on $N$ and $2N$ mesh sizes respectively. And the order of convergence is calculated by the formula $p_{i}=\log _{2} \left(\frac{D_{i}^{N}}{D_{i}^{2N}}\right)$.

\begin{figure} 
\centering
\includegraphics[scale=0.5]{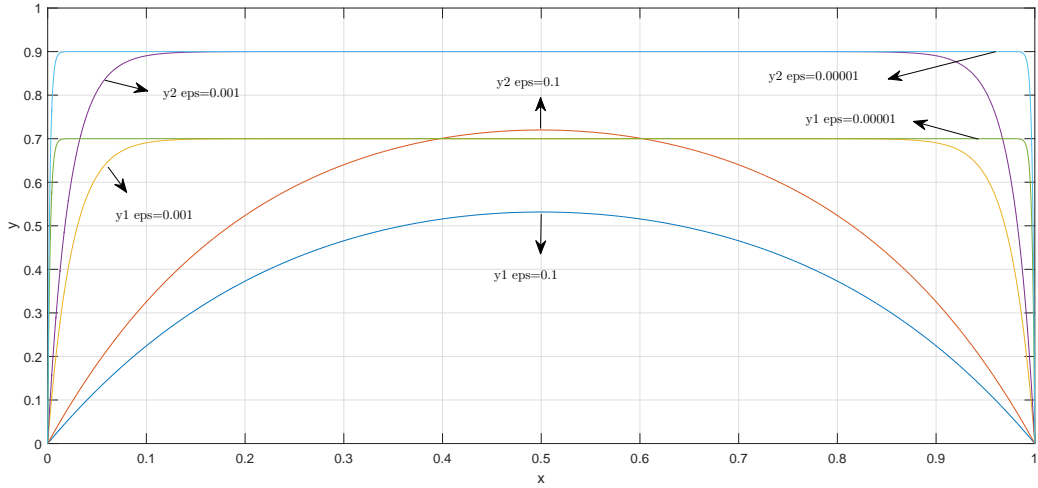}
\caption{\label{figure1_ex1} SCEM approximations for
illustrative Example \ref{exam1} for certain values of $\protect\varepsilon $.}
\label{fig:1}
\end{figure}

\begin{table} 
\centering
\begin{tabular}{|c|c|c|c|}
\hline
$x$ & $y_{1}^{\text{hybrid}}$ for$\ \varepsilon =1$ & $y_{1}^{\text{hybrid}}$
for$\ \varepsilon =0.01$ & $y_{1}^{\text{hybrid}}$ for$\ \varepsilon =0.0001$ \\
\hline
$0.000$ & $0.000000000000000$ & $0.000000000000000$ & $0.000000000000000$ \\
\hline
$0.001$ & $0.000476481355558$ & $0.008753818527279$ & $0.083181084554497$ \\
\hline
$0.003$ & $0.001426444995724$ & $0.025964925472142$ & $0.222962651683962$ \\
\hline
$0.070$ & $0.030951927479662$ & $0.418207225230131$ & $0.699958187708815$ \\
\hline
$0.090$ & $0.038906391001189$ & $0.484445481559528$ & $0.699997527634950$ \\
\hline
$0.100$ & $0.042736593719825$ & $0.511650798815475$ & $0.699999398893554$ \\
\hline
$0.300$ & $0.099085092371664$ & $0.688146304878580$ & $0.700000000000000$ \\
\hline
$0.500$ & $0.117696173594857$ & $0.698588175505725$ & $0.700000000000000$ \\
\hline
$0.700$ & $0.099085092371662$ & $0.688146304878580$ & $0.700000000000000$ \\
\hline
$0.900$ & $0.042736593719823$ & $0.511650798815475$ & $0.699999398893554$ \\
\hline
$0.910$ & $0.038906391001187$ & $0.484445481559528$ & $0.699997527634950$ \\
\hline
$0.930$ & $0.030951927479661$ & $0.418207225230131$ & $0.699958187708815$ \\
\hline
$0.997$ & $0.001426444995724$ & $0.025964925472143$ & $0.222962651683962$ \\
\hline
$0.999$ & $0.000476481355559$ & $0.008753818527279$ & $0.083181084554497$ \\
\hline
$1.000$ & $0.000000000000000$ & $0.000000000000000$ & $0.000000000000000$ \\
\hline
\end{tabular}
\caption{\label{table1} Approximations to $y_{1}$ of illustrative Example \ref{exam1} for various values of $ \varepsilon $, and $N=1024$.}
\end{table}

\begin{table} 
\centering
\begin{tabular}{|c|c|c|c|}
\hline
$x$ & $y_{2}^{\text{hybrid}}$ for $\varepsilon =1$ & $y_{2}^{\text{hybrid}}$
for $\varepsilon =0.01$ & $y_{2}^{\text{hybrid}}$ for $\varepsilon =0.0001$ \\
\hline
$0.000$ & $0.000000000000000$ & $0.000000000000000$ & $0.000000000000000$ \\
\hline
$0.001$ & $0.000836831145262$ & $0.013176322987849$ & $0.123255082602135$ \\
\hline
$0.003$ & $0.002504501565593$ & $0.038941217642193$ & $0.320704685367825$ \\
\hline
$0.070$ & $0.053859654789803$ & $0.576399973393426$ & $0.899958155847379$ \\
\hline
$0.090$ & $0.067538818218910$ & $0.657713587952894$ & $0.899997527270993$ \\
\hline
$0.100$ & $0.074102135126685$ & $0.690275199692044$ & $0.899999398854655$ \\
\hline
$0.300$ & $0.168930685338889$ & $0.887902110952330$ & $0.900000000000000$ \\
\hline
$0.500$ & $0.199554186019598$ & $0.898582598753880$ & $0.900000000000000$ \\
\hline
$0.700$ & $0.168930685338888$ & $0.887902110952330$ & $0.900000000000000$ \\
\hline
$0.900$ & $0.074102135126684$ & $0.690275199692044$ & $0.899999398854655$ \\
\hline
$0.910$ & $0.067538818218910$ & $0.657713587952894$ & $0.899997527270993$ \\
\hline
$0.930$ & $0.053859654789803$ & $0.576399973393426$ & $0.899958155847379$ \\
\hline
$0.997$ & $0.002504501565593$ & $0.038941217642193$ & $0.320704685367825$ \\
\hline
$0.999$ & $0.000836831145262$ & $0.013176322987849$ & $0.123255082602135$ \\
\hline
$1.000$ & $0.000000000000000$ & $0.000000000000000$ & $0.000000000000000$ \\
\hline
\end{tabular}
\caption{\label{table2} Approximations to $y_{2}$ of Example \ref{exam1} for various values of $\varepsilon $ and $N=1024$.}
\end{table}

\begin{table} 
\centering
\begin{tabular}{|c|c|c|c|c|c|}
\hline
$\varepsilon $ & $N=64$ & $N=128$ & $N=256$ & $N=512$ & $N=1024$ \\ \hline
$2^{-1}$ & $6.12679e-10$ & $3.70626e-11$ & $2.27984e-12$ & $1.40748e-13$ & $%
6.27276e-15$ \\ \hline
$2^{-2}$ & $1.81544e-09$ & $1.09791e-10$ & $6.76309e-12$ & $4.16611e-13$ & $%
3.53050e-14$ \\ \hline
$2^{-3}$ & $6.22012e-09$ & $3.76598e-10$ & $2.31717e-11$ & $1.43662e-12$ & $%
8.91509e-14$ \\ \hline
$2^{-4}$ & $2.47671e-08$ & $1.50369e-09$ & $9.25755e-11$ & $5.73113e-12$ & $%
3.58879e-13$ \\ \hline
$2^{-5}$ & $1.01353e-07$ & $6.14983e-09$ & $3.78493e-10$ & $2.34682e-11$ & $%
1.45650e-12$ \\ \hline
$2^{-6}$ & $4.03660e-07$ & $2.46994e-08$ & $1.52268e-09$ & $9.43695e-11$ & $%
5.87935e-12$ \\ \hline
$2^{-7}$ & $1.61007e-06$ & $9.89155e-08$ & $6.09315e-09$ & $3.77304e-10$ & $%
2.34998e-11$ \\ \hline
$2^{-8}$ & $6.44894e-06$ & $3.92055e-07$ & $2.43357e-08$ & $1.51185e-09$ & $%
9.41339e-11$ \\ \hline
$2^{-9}$ & $2.43061e-05$ & $1.56232e-06$ & $9.74256e-08$ & $6.04717e-09$ & $%
3.76837e-10$ \\ \hline
$2^{-10}$ & $4.88247e-07$ & $6.24635e-06$ & $3.86149e-07$ & $2.41664e-08$ & $%
1.50718e-09$ \\ \hline
$2^{-11}$ & $1.42363e-06$ & $2.36149e-05$ & $1.53878e-06$ & $9.66679e-08$ & $%
6.02468e-09$ \\ \hline
$2^{-12}$ & $4.88247e-07$ & $2.39459e-07$ & $6.14740e-06$ & $3.83279e-07$ & $%
2.40598e-08$ \\ \hline
$2^{-13}$ & $2.20197e-05$ & $7.04572e-07$ & $2.32768e-05$ & $1.52714e-06$ & $%
9.62998e-08$ \\ \hline
$2^{-14}$ & $1.92995e-05$ & $2.39028e-07$ & $1.18093e-07$ & $6.09883e-06$ & $%
3.81771e-07$ \\ \hline
$2^{-15}$ & $4.12438e-07$ & $2.16196e-05$ & $3.43788e-07$ & $2.31095e-05$ & $%
1.52132e-06$ \\ \hline
$D_{1}^{N}$ & $2.43061e-05$ & $2.36149e-05$ & $2.32768e-05$ & $2.31095e-05$
& $1.52132e-06$ \\ \hline
$p_{1}^{N}$ & $0.0416210001$ & $0.0208046902$ & $0.0312113804$ & $%
3.9250904939$ & $3.9889941035$ \\ \hline
\end{tabular}
\caption{ \label{table3} $D_{1}^{N}$ and $p_{1}^{N}$ for illustrative Example \ref{exam1} for various values of $\varepsilon $ and $N$.}
\end{table}

\begin{figure} 
\centering
\includegraphics[scale=0.5]{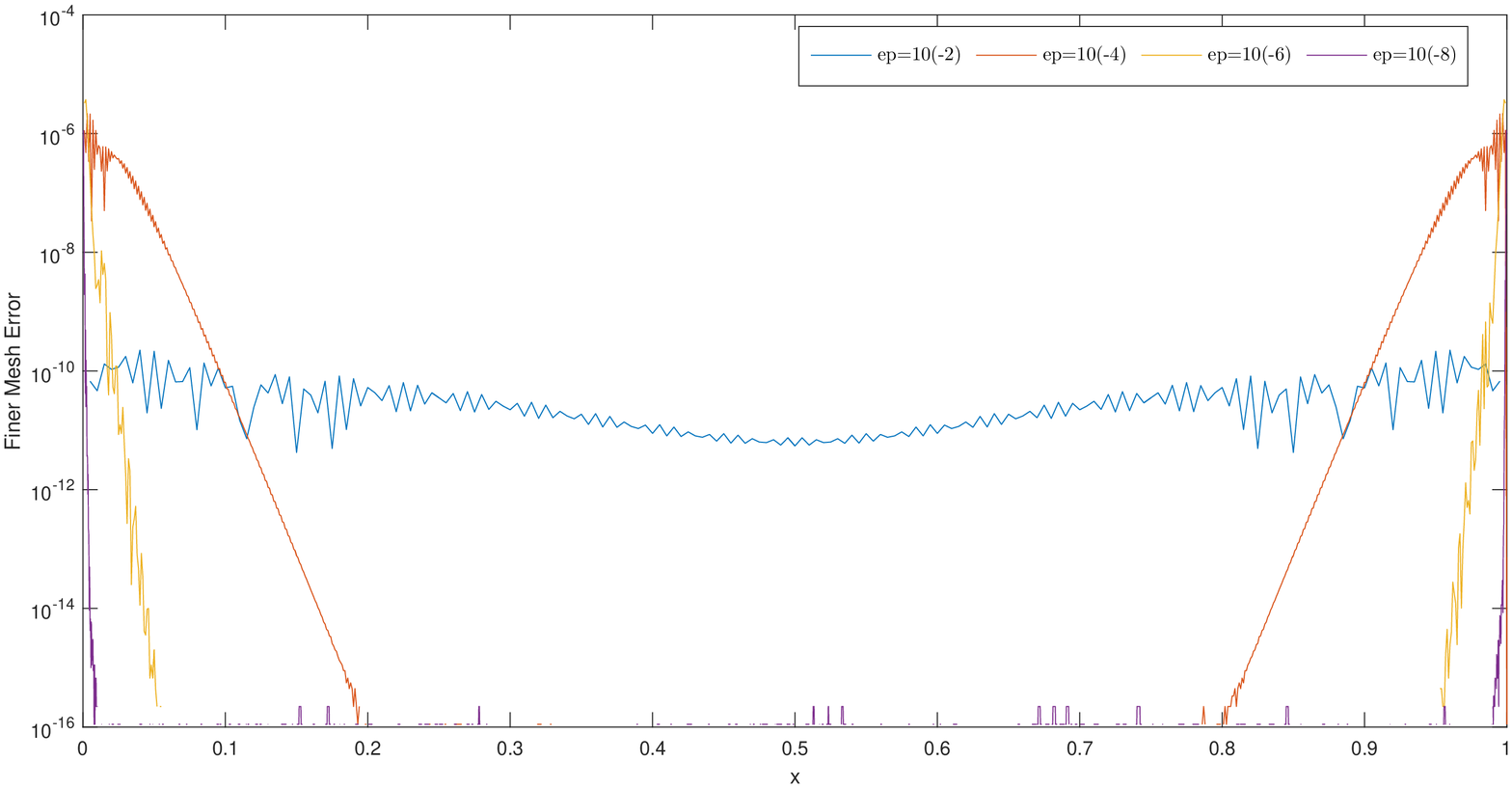}
\caption{\label{figure2_ex1} Errors in $y_{1}$ approximations of Example \ref{exam1} for various values of $\protect\varepsilon $.}
\label{fig:2}
\end{figure}

\begin{example}\label{exam2}
Consider the system of coupled singularly perturbed reaction-diffusion equations \cite{17}:
\begin{equation}\label{exam_pb2}
\left\{
\begin{array}{ll}
-\varepsilon y_{1}''\left( x\right) +3y_{1}\left( x\right)
- y_{2}\left( x\right) -y_{3}\left( x\right) = 0, \\[8pt]
-\varepsilon y_{2}''\left( x\right) -y_{1}\left( x\right)
+ 3y_{2}\left( x\right) -y_{3}\left( x\right) = 1, \\[8pt]
-\varepsilon y_{3}''\left( x\right) -y_{2}\left( x\right)
+ 3y_{3}\left( x\right) = x, \\[8pt]
y_{1}\left( 0\right) = y_{1}\left(1\right) = 0, \quad y_{2}\left( 0\right) =y_{2}\left( 1\right) =0, \quad y_{3}\left( 0\right) =y_{3}\left( 1\right) =0.
\end{array}\right.
\end{equation}
\end{example}
Solution of this problem exhibits boundary layer behavior in all the
components $y_{1},$ $y_{2}$ and $y_{3}.$

\begin{figure} 
\centering
\includegraphics[scale=0.5]{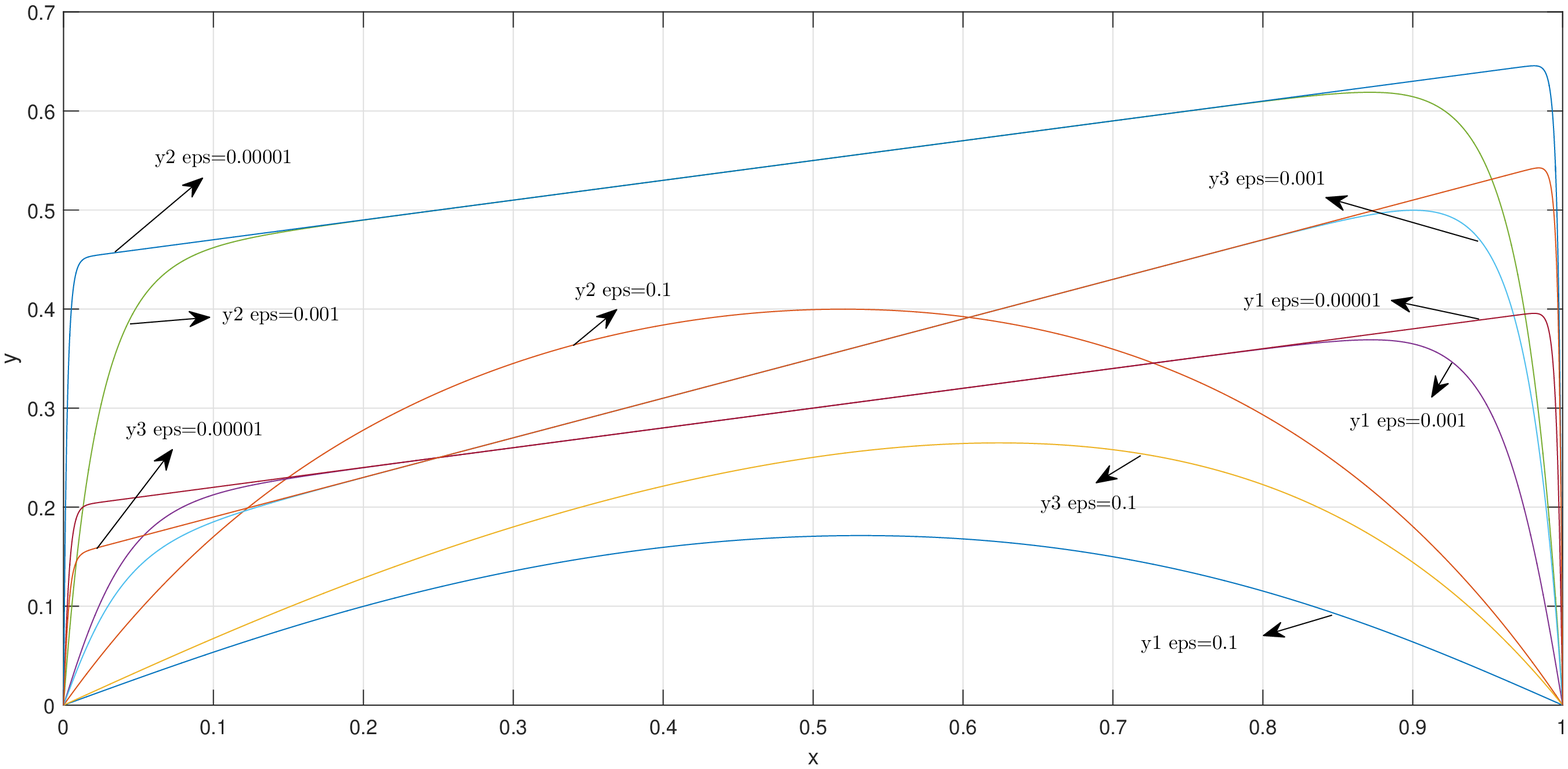}
\caption{\label{figure1_ex2} SCEM approximations for Example \ref{exam2} for certain values of $\protect\varepsilon$.}
\label{fig:3}
\end{figure}

\begin{table} 
\centering
\begin{tabular}{|c|c|c|c|}
\hline
$x$ & $y_{3}^{\text{hybrid}}$ for $\varepsilon =1$ & $y_{3}^{\text{hybrid}}$
for $\varepsilon =0.01$ & $y_{3}^{\text{hybrid}}$ for $\varepsilon =0.0001$ \\
\hline
$0.000$ & $0.000000000000000$ & $0.000000000000000$ & $0.000000000000000$ \\
\hline
$0.001$ & $0.000296931914224$ & $0.001555243868093$ & $0.011907872130351$ \\
\hline
$0.003$ & $0.000445395437314$ & $0.004664568760707$ & $0.034807883220829$ \\
\hline
$0.070$ & $0.010337895722462$ & $0.098994889589599$ & $0.177946516870041$ \\
\hline
$0.090$ & $0.013246127387130$ & $0.121783382519711$ & $0.185994913258379$ \\
\hline
$0.100$ & $0.014688076281396$ & $0.132245460929406$ & $0.189998430525371$ \\
\hline
$0.300$ & $0.040416458359500$ & $0.263967353662464$ & $0.270000000000000$ \\
\hline
$0.500$ & $0.055600510502286$ & $0.348297797967831$ & $0.350000000000000$ \\
\hline
$0.700$ & $0.053365373735929$ & $0.417475745785387$ & $0.430000000000000$ \\
\hline
$0.900$ & $0.025982747779700$ & $0.360670915745580$ & $0.509996860535433$ \\
\hline
$0.910$ & $0.023793671058585$ & $0.344253981906367$ & $0.513989822709124$ \\
\hline
$0.930$ & $0.019152362750862$ & $0.302349198379968$ & $0.521892825852710$ \\
\hline
$0.997$ & $0.000917731428012$ & $0.020288336492293$ & $0.178812725352010$ \\
\hline
$0.999$ & $0.000306908478447$ & $0.006860930176790$ & $0.067932980344570$ \\
\hline
$1.000$ & $0.000000000000000$ & $0.000000000000000$ & $0.000000000000000$ \\
\hline
\end{tabular}
\caption{\label{table4} Approximations to $y_{3}$ of Example \ref{exam2} for various values of $\varepsilon $, and $N=1024$.}
\end{table}

\begin{table} 
\centering
\begin{tabular}{|c|c|c|c|c|c|}
\hline
$\varepsilon $ & $N=64$ & $N=128$ & $N=256$ & $N=512$ & $N=1024$ \\ \hline
$2^{-1}$ & $3.97376e-08$ & $2.40516e-09$ & $1.48046e-10$ & $9.04537e-12$ & $%
5.14129e-13$ \\ \hline
$2^{-2}$ & $1.60003e-07$ & $9.68587e-09$ & $5.95688e-10$ & $3.69476e-11$ & $%
2.45948e-12$ \\ \hline
$2^{-3}$ & $4.55184e-07$ & $2.75605e-08$ & $1.69540e-09$ & $1.05172e-10$ & $%
6.50464e-12$ \\ \hline
$2^{-4}$ & $1.86005e-06$ & $1.12657e-07$ & $6.93077e-09$ & $4.29730e-10$ & $%
2.67455e-11$ \\ \hline
$2^{-5}$ & $5.30766e-06$ & $3.21668e-07$ & $1.97922e-08$ & $1.22731e-09$ & $%
7.63886e-11$ \\ \hline
$2^{-6}$ & $7.19988e-05$ & $1.28455e-06$ & $7.90643e-08$ & $4.90311e-09$ & $%
3.05273e-10$ \\ \hline
$2^{-7}$ & $8.50572e-05$ & $3.62682e-06$ & $2.23376e-07$ & $1.38547e-08$ & $%
8.62547e-10$ \\ \hline
$2^{-8}$ & $1.48198e-04$ & $4.15936e-05$ & $8.93032e-07$ & $5.54076e-08$ & $%
3.44981e-09$ \\ \hline
$2^{-9}$ & $1.20331e-04$ & $4.88494e-05$ & $2.52367e-06$ & $1.56680e-07$ & $%
9.75681e-09$ \\ \hline
$2^{-10}$ & $1.12466e-05$ & $1.02047e-04$ & $2.89999e-05$ & $6.26449e-07$ & $%
3.90227e-08$ \\ \hline
$2^{-11}$ & $1.42610e-05$ & $8.24795e-05$ & $3.41249e-05$ & $1.77034e-06$ & $%
1.10348e-07$ \\ \hline
$2^{-12}$ & $1.25522e-04$ & $3.92817e-06$ & $7.24787e-05$ & $1.73510e-05$ & $%
4.41207e-07$ \\ \hline
$2^{-13}$ & $7.95760e-05$ & $4.97760e-06$ & $7.11641e-05$ & $3.54793e-05$ & $%
1.24685e-06$ \\ \hline
$2^{-14}$ & $5.17968e-05$ & $8.69193e-05$ & $1.41263e-06$ & $5.00248e-05$ & $%
1.43538e-05$ \\ \hline
$2^{-15}$ & $5.35070e-06$ & $5.53248e-05$ & $1.76076e-06$ & $5.00871e-05$ & $%
2.50118e-05$ \\ \hline
$D_{3}^{N}$ & $1.48198e-04$ & $1.02047e-04$ & $7.24787e-05$ & $5.00871e-05$
& $2.50118e-05$ \\ \hline
$p_{3}^{N}$ & $0.5382881157$ & $0.4936047862$ & $0.5331179957$ & $%
1.0018301969$ & $4.7990910350$ \\ \hline
\end{tabular}
\caption{ \label{table5} $D_{3}^{N}$ and $p_{3}^{N}$ for Example \ref{exam2}, for various values of $\varepsilon $ and $N$.}
\end{table}

\begin{figure}
\centering
\includegraphics[scale=0.5]{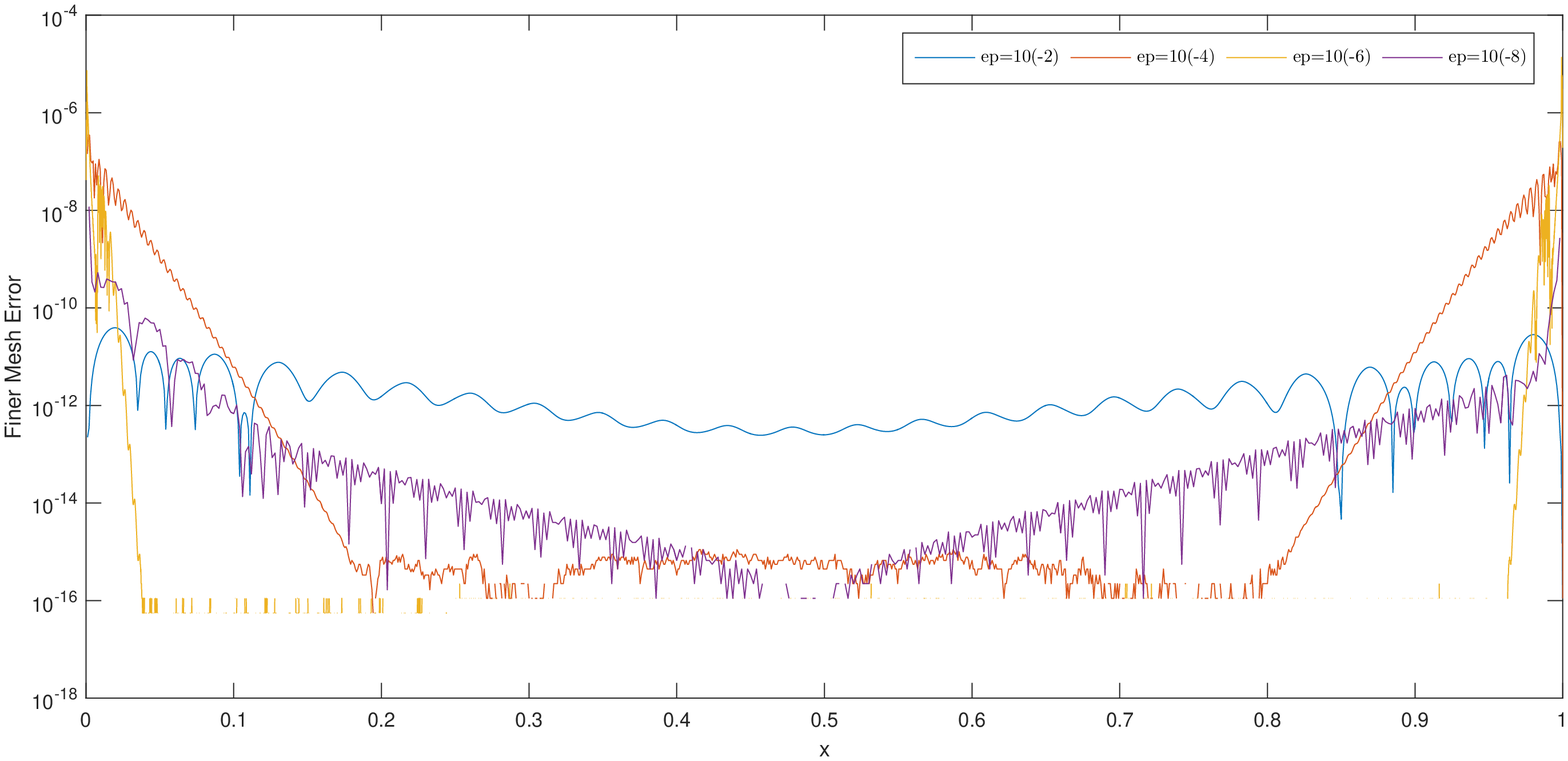}
\caption{\label{fig4} Errors in $y_{3}$ approximations of Example \ref{exam2} for various values of $\protect\varepsilon $.}
\label{fig:mesh}
\end{figure}

\section{Conclusion}

In this paper, singularly perturbed  system of two-point reaction-diffusion type boundary-value problems are examined. In order to obtain good approximations to the solution of these types of problems, a hybrid method which consists of an asymptotic method, known as SCEM and a numerical method based on finite differences given in \cite{25} is proposed. In Section \ref{sec4}, the implementation of the present method is given in detail on the Illustrative Example \ref{exam1}. In Table \ref{table1}, Table \ref{table2} and Table \ref{table4} numerical results that obtained by the hybrid method and in Table \ref{table3} and Table \ref{table5} maximum point-wise errors and convergence based on finer mesh strategy are presented. Additionally, approximations for some $\varepsilon $ values are given in Figure \ref{figure1_ex1} and Figure \ref{figure1_ex2}  to illustrate the layer behavior. Maximum point-wise errors for some $\varepsilon $ values are given in Figure \ref{figure2_ex1} and Figure \ref{fig4}. As one can point out from the figures and tables, that the hybrid method gives high-accurate results and well-suited for singularly perturbed system of reaction-diffusion equations.

\end{document}